\documentclass{amsart}

\usepackage{amssymb}
\usepackage{amsthm}
\usepackage[matrix,arrow,graph,curve,frame]{xy}
\usepackage{hyperref}
\usepackage{mathrsfs}

\theoremstyle{definition} 
\theoremstyle{definition} \newtheorem*{lemma}{Lemma}
\theoremstyle{definition} 
\theoremstyle{definition} \newtheorem*{example}{Example}

\newcommand{\A}{\ensuremath{\mathscr{A}}}
\newcommand{\C}{\ensuremath{\mathscr{C}}}
\newcommand{\V}{\ensuremath{\mathscr{V}}}
\newcommand{\ox}{\ensuremath{\otimes}}
\newcommand{\op}{\ensuremath{\mathrm{op}}}
\newcommand{\ob}{\ensuremath{\mathrm{ob}}}
\newcommand{\Vect}{\ensuremath{\mathbf{Vect}}}

\newcommand{\ra}{\ensuremath{\xymatrix@C=4.5ex@1{\ar[r]&}}}
\newcommand{\Ra}{\ensuremath{\xymatrix@C=5ex@1{\ar@{=>}[r]&}}}

\title{Compact convolution}
\author{Brian J. Day}
\address{Centre of Australian Category Theory, Macquarie University, NSW,
2109, Australia}
\date{March 27, 2006}

\begin{document}

\begin{abstract}
We state a Yoneda-type lemma which leads to various functor categories being
compact closed.
\end{abstract}

\maketitle

\begin{lemma}
Given a $\V$-functor
\[
T:\A^\op \ox \A \ra \V
\]
with $\A$ a $\V_f$-category, suppose that the canonical map
\[
\int^X \int_Y \A(Y,X) \ox T(X,Y) \ra \int_Y \int^X \A(Y,X) \ox T(X,Y)
\]
is an isomorphism. Then, for each choice of $\V$-natural isomorphism
\[
\A(Y,X) \cong \A(X,Y)^*,
\]
we get an isomorphism
\[
\int^X T(X,X) \xymatrix@1{\ar[r]_{\alpha}^{\cong} &} \int_Y T(Y,Y)
\]
where $\alpha_{XY}$ is the $\V$-dinatural composite
\[
\xymatrix{T(X,X) \ar[r]^-{\mathrm{can}} & [\A(X,Y),T(X,Y)] \ar[r]^-{\cong} &
\A(Y,X) \ox T(X,Y) \ar[r]^-{\mathrm{can}} & T(Y,Y)}.
\]
\end{lemma}

\begin{example}
Let $\V = \Vect_k$ and suppose $\ob(\A)$ is finite; then $\int^X S(X,X)$ is
absolute for all $\V$-functors $S:\A^\op \ox \A \ra \V$, hence
\[
\int^X \int_Y \xymatrix@1{\ar[r]^{\cong} &} \int_Y \int^X
\]
always, so
\[
\int^X T(X,X) \xymatrix@1{\ar[r]^{\cong} &} \int_Y T(Y,Y)
\]
for $\A = k_*(\text{finite groupoid})$ or $\A$ a finite dimensional Hopf
algebra, etc.
\end{example}

\subsubsection*{Extension}

Suppose $\V = \Vect_k$ and $\A \subset \C$ is $\V$-dense with $\ob(\A)$
finite; then, for suitably continuous $\V$-functors
\[
T:\C^\op \ox \C \ra \V
\]
we have
\[
\int^C T(C,C) \xymatrix@1{\ar[r]^{\cong} &} \int_D T(D,D)
\]
for each choice of $\V$-natural isomorphism $\A(X,Y) \cong \A(Y,X)^*$.

\begin{proof}
Both
\[
\int^X T(X,X) \xymatrix@1{\ar[r]^{\cong} &} \int^C T(C,C)
\quad \text{and} \quad
\int_D T(D,D) \xymatrix@1{\ar[r]^{\cong} &} \int_Y T(Y,Y)
\]
are isomorphisms by hypothesis and Yoneda.
\end{proof}

\subsubsection*{Convolution}

Let
\[
p:\A^\op \ox \A^\op \ox \A \ra \V_f \quad \text{and} \quad
j:\A \ra \V_f
\]
be a (commutative) promonoidal category where $p$, $\A(-,-)$, and $j$ have
finite support in each variable separately, and suppose there is a natural
isomorphism
\begin{equation}\label{*}\tag{$*$}
\int^{XY} j(Y) \ox p(X,B,Y)^* \ox p(X,C,A) \cong p(A,B,C)^*.
\end{equation}

For each $\V$-functor $G:\A \ra \V_f$ define $G^*:\A \ra \V_f$ by
\[
G^*(X) = \int^B (GB)^* \ox \int^Y j(Y) \ox p(X,B,Y)^*.
\]
Then
\begin{align*}
(G^* \ox H)(A)
  &= \int^{XC} G^*X \ox HC \ox p(X,C,A) \qquad \text{(by defn of $\ox$)}\\
  &= \int^{XC} \int^{B} (GB)^* \ox \int^Y j(Y) \ox p(X,B,Y)^* \ox HC
\ox p(X,C,A) \\
  &\cong \int^{BC} (GB)^* \ox HC \ox p(A,B,C)^* \qquad \text{(by \eqref{*})}\\
  &\cong \int_{BC} (GB \ox p(A,B,C))^* \ox HC \qquad \text{(by the lemma)}\\
  &\cong \int_{BC} [GB \ox p(A,B,C),HC] \\
  &= [G,H](A).
\end{align*}
So $[\A,\V_f]$ is compact.

\begin{example}
If $\A$ has an antipode
\[
S:\A^\op \ra \A,
\]
with $S^2 \cong 1$, then \eqref{*} holds if both
\[
p(X,Y,Z)^* \cong p(SX,SY,SZ) \quad \text{and} \quad p(X,Y,SZ) \cong p(Y,Z,SX).
\]

\begin{proof}
\begin{align*}
\int^{XY} & j(Y) \ox p(X,B,Y)^* \ox p(X,C,A) \\
&\cong \int^{XY} j(Y) \ox p(SX,SB,SY) \ox p(X,C,A) && \text{by hyp,}\\
&\cong \int^{XY} j(Y) \ox p(SB,Y,X) \ox p(X,C,A) && \text{by hyp and $S^2
\cong 1$,} \\
&\cong \int^{X} \A(SB,X) \ox p(X,C,A) && \text{since $j *p \cong \A(-,-)$,}\\
&\cong p(SB,C,A) && \text{by Yoneda,} \\
&\cong p(SA,SB,SC)  && \text{by hyp and $S^2 \cong 1$,} \\
&\cong p(A,B,C)^* && \text{by hyp.}
\qedhere
\end{align*}
\end{proof}
\end{example}

A special case of the above is where $\ob(\A)$ is finite. In the literature
(cf.~\cite{1} and~\cite{2}) the situation occurs where $\A \subset \C$ with
$(\C,\ox,I)$ a commutative compact closed category and $\A$ has the trace
promonoidal structure induced by $\C$ and $I \in \A$; that is
\[
p(X,Y,Z) = \C(X \ox Y,Z) \quad \text{and} \quad j(X) = \C(I,X).
\]
Of course we suppose that $\A(X,Y) = \A(Y,X)^*$, but we must also suppose
that 
\[
\C(X \ox Y,Z) \cong \C(Z,X \ox Y)^*
\]
naturally in $X,Y,Z \in \A$ in order for
\begin{align*}
p(SX,SY,SZ) &= \C(SX \ox SY,SZ) \\
            &= \C(Z, X \ox Y) \\
            &= \C(X \ox Y, Z)^* \\
            &= p(X,Y,Z)^*,
\end{align*}
and
\begin{align*}
p(X,Y,SZ) &= \C(X \ox Y,SZ) \\
          &= \C(X \ox Y \ox Z,I) \\
          &= \C(Y \ox Z, SX) \\
          &= p(Y,Z,SX).
\end{align*}

The point is that the empirical $\A$, being finite, tends not to be closed
under the $\ox$ of $\C$.

\begin{example}
Let $\V = \Vect_k$. Suppose $\A \subset \C$ is Cauchy dense; then
\[
[\C,\V_f] \simeq [\A,\V_f].
\]
If $\C$ is a commutative compact closed category and $\ob(\A)$ is finite
then $[\A,\V_f]$ is compact closed (as we saw earlier) so $[\C,\V_f]$ is
compact.
\end{example}



\end{document}